\documentclass{article}
\usepackage{cite}
\usepackage{amsmath,amssymb,amsfonts}
\usepackage{algorithmic}
\usepackage{graphicx}
\usepackage{textcomp}
\usepackage{xcolor}
\usepackage[justification=centering]{caption}
\usepackage[justification=centering]{subcaption}
\usepackage[font=footnotesize]{caption}
\usepackage{float}
\usepackage{pifont}
\usepackage{xcolor}
\definecolor{Alive}{rgb}{0.5, 0, 0.5}
\definecolor{Dead}{rgb}{1, 0.5, 0}
\definecolor{Links}{HTML}{7A0022}
\usepackage{algorithm2e}\RestyleAlgo{ruled}
\usepackage{algorithmic}
\usepackage{mathtools}
\usepackage{stackengine}
\usepackage{hyperref}
\hypersetup{
    colorlinks=true,
    citecolor=Links,
    linkcolor=Links,
    filecolor=Links,      
    urlcolor=Links,
    pdftitle={Organic Structures Emerging From Bio-Inspired Graph-Rewriting Automata},
    pdfpagemode=FullScreen,
}
  
\def\BibTeX{{\rm B\kern-.05em{\sc i\kern-.025em b}\kern-.08em
    T\kern-.1667em\lower.7ex\hbox{E}\kern-.125emX}}
    
\providecommand{\keywords}[1]
{
  \small	
  \textbf{\textit{Keywords ---}} #1
}
    
\begin{document}

\title{Organic Structures Emerging From Bio-Inspired Graph-Rewriting Automata}

\author{\href{https://orcid.org/0000-0002-3866-7615}{Paul Cousin} and  \href{https://orcid.org/0000-0002-0905-2515}{Aude Maignan}  \\
        \small Univ. Grenoble Alpes, CNRS, Grenoble INP, LJK \\
        \small Grenoble, France
}

\date{}

\maketitle

\vspace{2mm}

\begin{abstract}
\;Graph-Rewriting Automata (GRA) are an extension of Cellular Automata to a dynamic structure using local graph-rewriting rules. This work introduces linear algebra based tools that allow for a practical investigation of their behavior in deeply extended time scales. A natural subset of GRA is explored in different ways thereby demonstrating the benefits of this method. Some elements of the subset were discovered to create chaotic patterns of growth and others to generate organic-looking graph structures. These phenomena suggest a strong relevance of GRA in the modeling natural complex systems. The approach presented here can be easily adapted to a wide range of GRA beyond the chosen subset. \vspace{2mm}
\end{abstract} 


\vspace{1.5cm}

\section{Introduction}
\label{section:introduction}
The idea of extending the study of discrete dynamical systems like cellular automata to systems with evolving typologies has already been approached in several ways. DEM-Systems \cite{edwards2020class,edwards2016systems,edwards2014complex}, Structurally Dynamic Cellular Automata \cite{ilachinski2009structurally,alonso2006structurally,alonso2007structurally}, Generative Network Automata \cite{sayama2007generative,sayama2009generative,schmidt2013designing} and Graph-Rewriting Automata \cite{tomita2007self, tomita2007asynchronous, tomita2009graph} are related concepts sharing this goal. These systems can produce a range of new behaviors compared to their fixed topology counterparts, however their implementation is much less straightforward. \\

In this paper, new linear algebra based tools will be presented, then applied to study a natural subset of GRA. Links to the Mathematica source code, a more \linebreak recent GPU-accelerated Python implementation, and more information can be found at \href{https://paulcousin.github.io/graph-rewriting-automata}{paulcousin.github.io/graph-rewriting-automata}.

\pagebreak

A GRA consists of an initial graph $G_0$ defined at time step $t=0$ and a rule to iteratively evolve it to any discrete time step $t>0$. $G_t$ is the graph obtained at time $t$. \linebreak The studied subset of GRA will be defined by the following restrictions:
\begin{enumerate}
    \item 3-regular, undirected and finite graphs,
    \item binary labeling of vertices $v$ with two possible states $s$:
    \begin{enumerate}
        \item[•] $s(v)=1$, called ``alive" and colored purple \,\parbox{2mm}{\includegraphics[width=2.6mm]{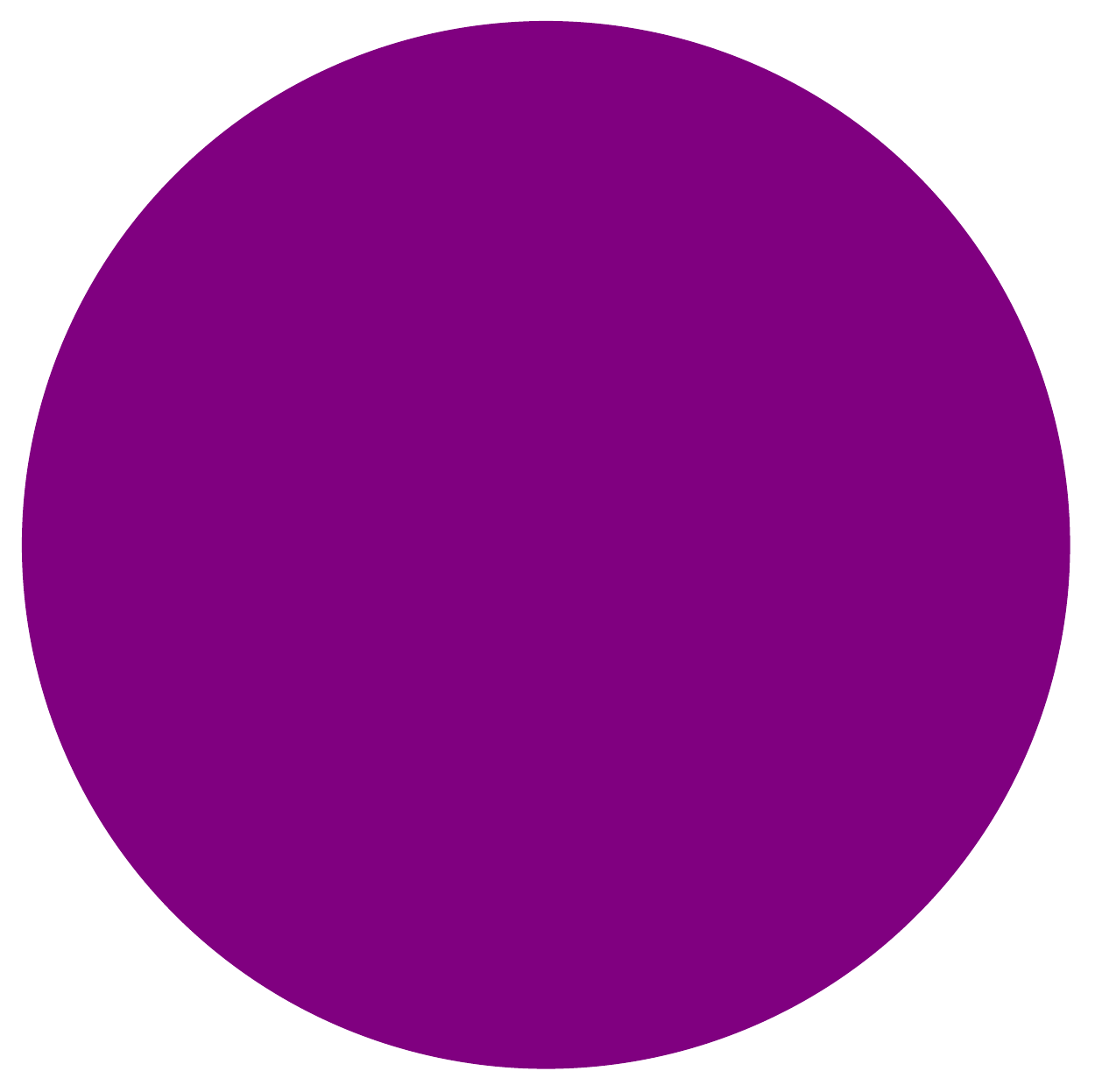}}\;\,,
        \item[•]  $s(v)=0$, called ``dead" and colored orange \,\parbox{2mm}{\includegraphics[width=2.6mm]{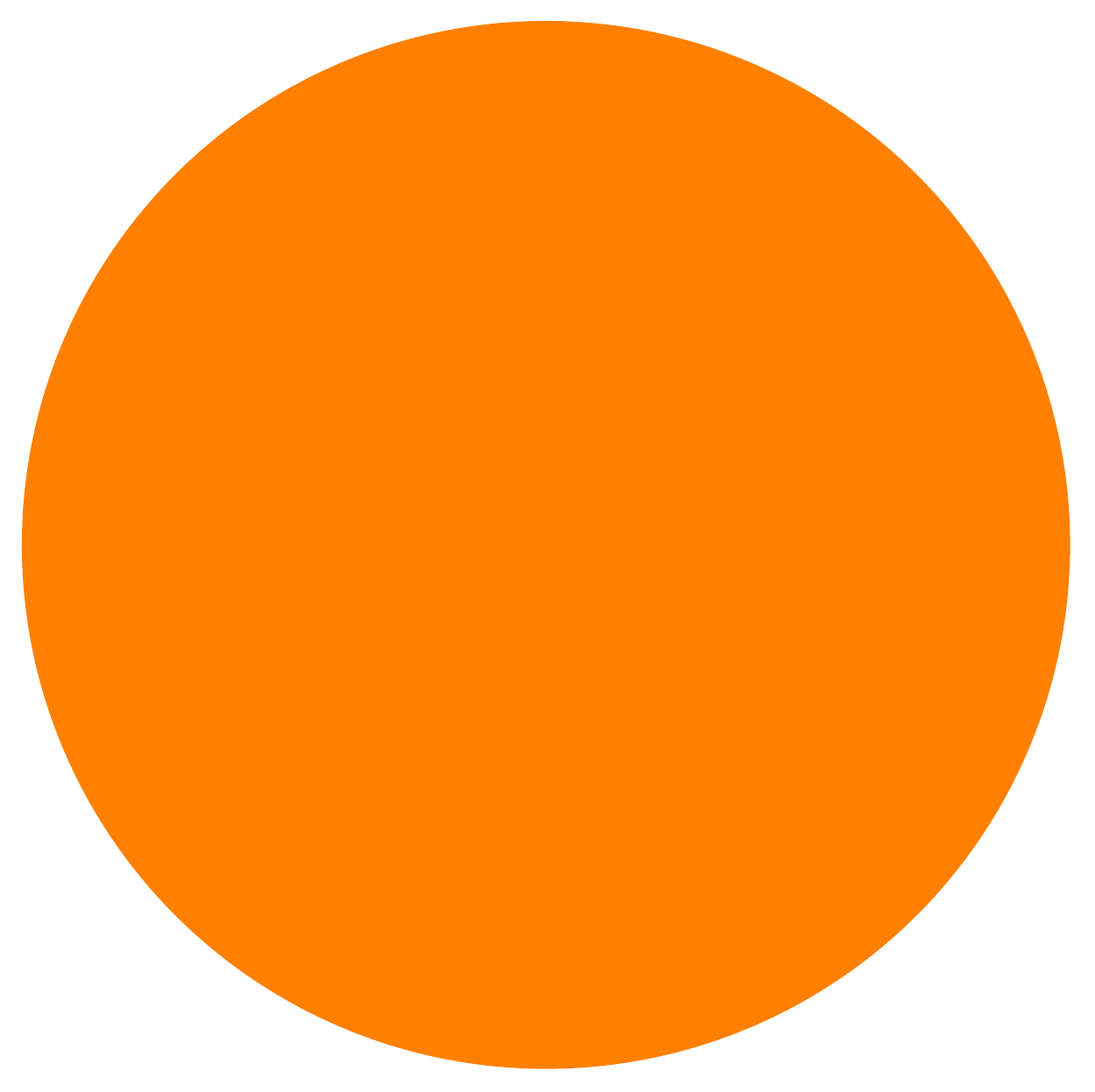}}\;\,,
    \end{enumerate}
    \item rules local to a vertex and its adjacent vertices, deterministic and applied simultaneously on the entire graph,
    \item only one type of topology altering operation called division.
\end{enumerate}

\noindent Divisions will be operated as follows: \parbox{2.5cm}{\includegraphics[height=1.2cm]{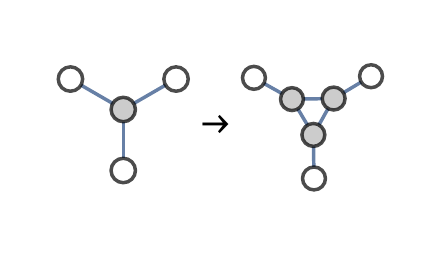}}

In this model, vertices can, for example, be seen as cells that divide under specific local conditions (internal and environmental). An initially simple graph can thus grow into a large organism with a complex self-organized structure.

\begin{figure}[H]
\begin{center}
\begin{tabular}{c c c c c c c c c c c } 
 \parbox{12mm}{\includegraphics[height=14mm,angle=180]{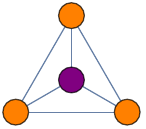}} &$\rightarrow$
 \parbox{12mm}{\includegraphics[height=14mm,angle=180]{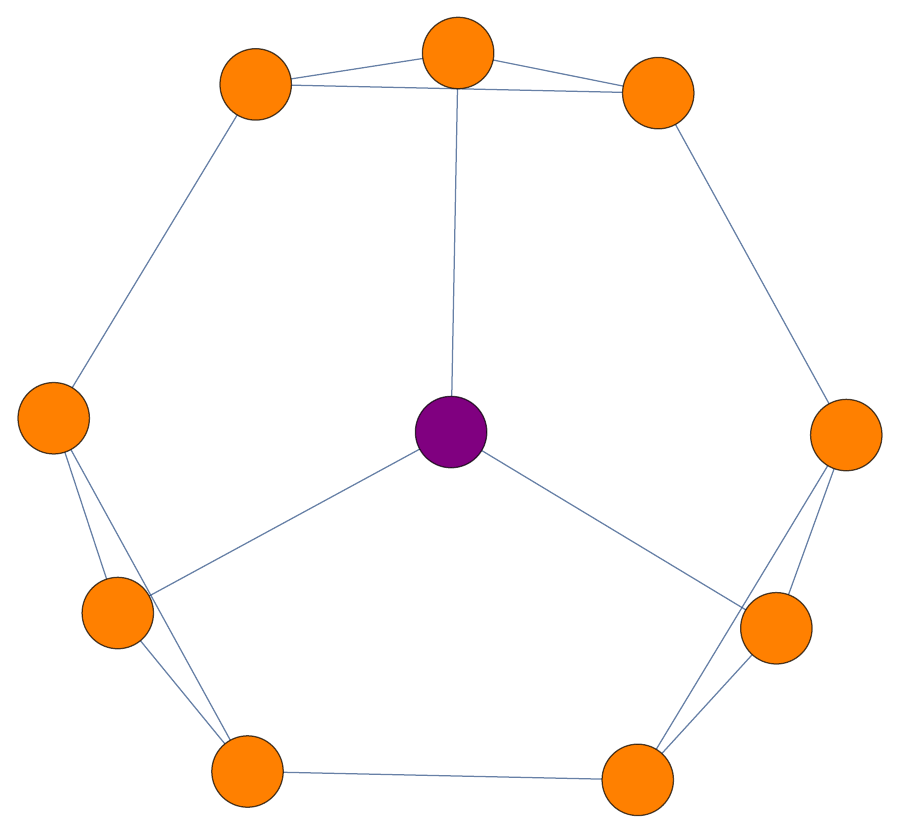}} &$\rightarrow$
 \parbox{12mm}{\includegraphics[height=14mm]{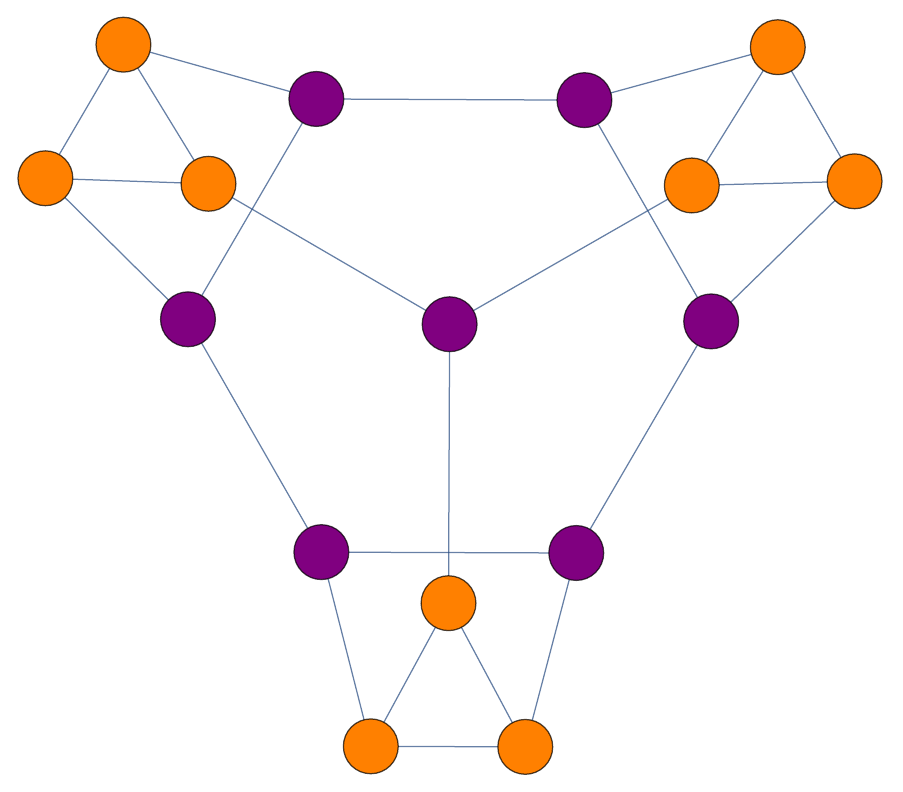}} &$\rightarrow$
 \parbox{12mm}{\includegraphics[height=14mm,angle=195]{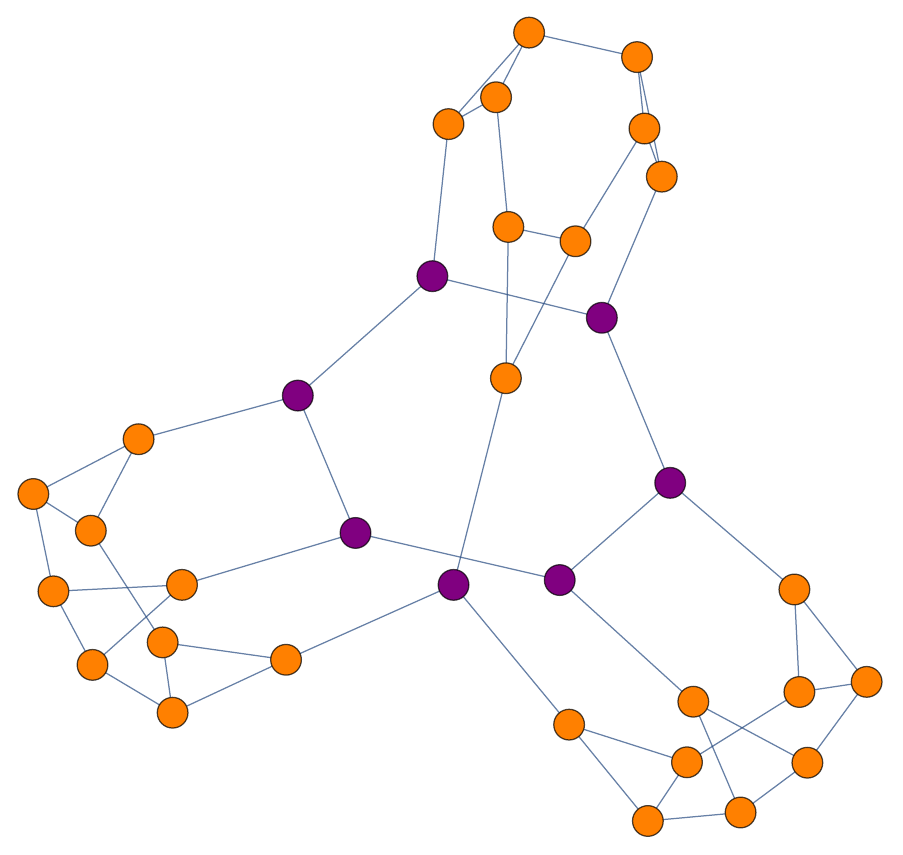}} &$\rightarrow$
 \parbox{12mm}{\includegraphics[height=14mm]{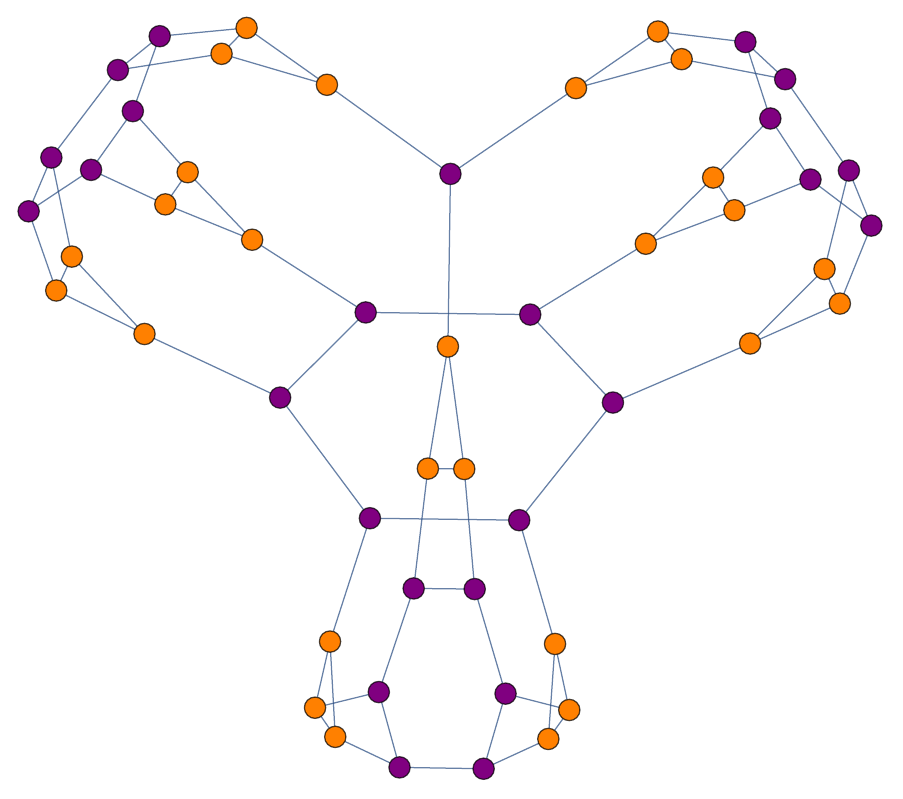}} & $\rightarrow\;\dots$ & \\
 \end{tabular}
 {\scriptsize  $G_0$
 \;\;\qquad\qquad\quad $G_1$
 \,\;\;\,\,\qquad\qquad\quad $G_2$
 \;\;\;\,\qquad\qquad\quad $G_3$
 \,\,\,\,\qquad\qquad\quad $G_4$
 \;\;\;\;\;\;\qquad\quad\quad}
\end{center}
\caption{Example of evolution starting from a simple initial graph.}
\label{example}
\end{figure}

\section{Graphs and Rules}

A \textbf{graph} $G$ is fully described by an \textbf{adjacency matrix} $\mathcal{A}$ and a \textbf{state vector} $\mathcal{S}$: \; $G_t=\{\,\mathcal{A}_t\, ,\,\mathcal{S}_t\,\}$

\begin{figure}[H]
\centering
\resizebox{8cm}{!}{
\begin{tabular}{c c c c c}
$G_0\equiv$\parbox{1.5cm}{\includegraphics[height=2cm]{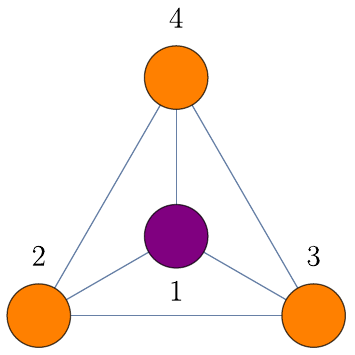}}
& \qquad\qquad & 
$\mathcal{A}_0=
\begin{pmatrix}
0 & 1 & 1 & 1 \\ 
1 & 0 & 1 & 1 \\ 
1 & 1 & 0 & 1 \\ 
1 & 1 & 1 & 0
\end{pmatrix}$ 
& &
$\mathcal{S}_0=
\begin{pmatrix}
1 \\ 0 \\ 0 \\ 0
\end{pmatrix}$ 
\end{tabular}}
\caption{Example of the graph $G_0$ from  \protect\textit{Figure \protect\ref{example}}.}
\label{example-graph}
\end{figure}

The \textbf{neighborhood} $N(v)$ of a \textbf{vertex} $v$ is the set of vertices which are adjacent to it. The \textbf{state} $s$ of a vertex being either dead or alive and there being from 0 to 3 alive vertices in its neighborhood, there are only 8 possible \textbf{configurations} $c$, which can be ordered in this way: 
\begin{equation}
c(v)=4\times s(v)+\sum_{i\in N(v)} s(i)
\end{equation}

\begin{figure}[H]
\centering
\includegraphics[width=9cm]{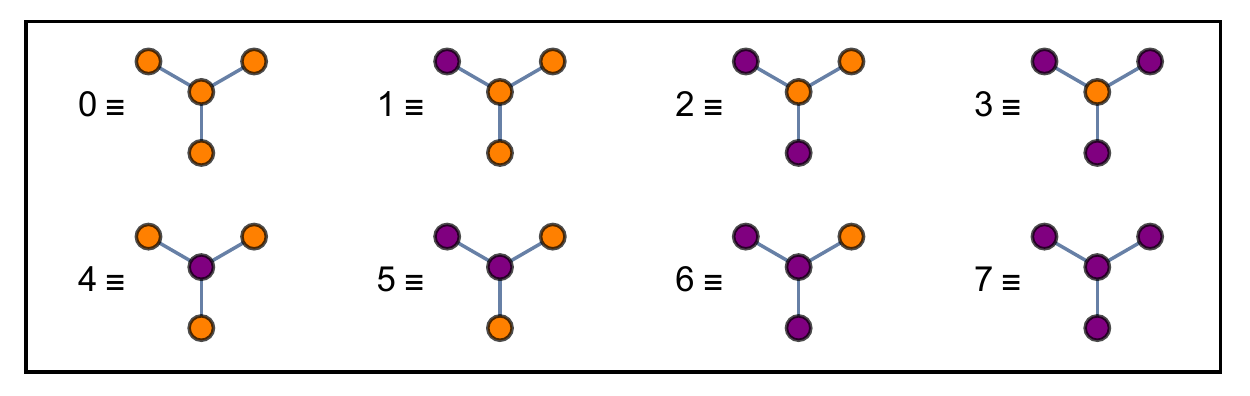}
\caption{Eight possible configurations.}
\label{configurations}
\end{figure}

For instance, if $v$ is a dead vertex surrounded by 3 dead vertices, $c(v)$ will be equal to 0 (see \textit{Figure \ref{configurations}} for all the other configurations).

The space of possible rules applicable in the subset of GRA defined in \textit{Section \ref{section:introduction}} is finite. Every rule must specify, for each configuration, whether the vertex will be alive or dead at $t+1$ and whether it will have undergone a division, leading to 4 possible final states. Thus, there are $4^8=65,536$ possible rules. Each rule can be described by two functions:
\begin{equation}
\begin{aligned}
& R: [[0,7]]\rightarrow\{\;0\equiv
\parbox{2mm}{\includegraphics[width=2.6mm]{graphics/rules/orange.pdf}}\;\,
,\;\, 1\equiv
\parbox{2mm}{\includegraphics[width=2.6mm]{graphics/rules/purple.pdf}}\;\,
\} 
\qquad\qquad\qquad\qquad\quad\; \\  
& R(c(v_t))=s(v_{t+1}) 
\end{aligned}
\end{equation}

\begin{equation}
\begin{aligned}
& R': [[0,7]]\rightarrow\{\;0\equiv\parbox{6mm}{\includegraphics[width=6mm]{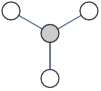}}\;,\;1\equiv\parbox{6mm}{\includegraphics[width=6mm]{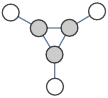}}\;\} \\  
& R'(c(v_t))=\begin{cases}\begin{array}{ll}
        1, & \text{if c leads to a division at t+1}\\
        0, & \text{otherwise}
        \end{array}\end{cases}
\end{aligned}
\end{equation}

\vspace{4mm}

\noindent Every rule can thus be labeled by a unique \textbf{rule number} $n$.
\begin{equation}
    n=\sum_{i=0}^7 \left[ 2^i R(i) +  2^{i+8} R'(i) \right]
\end{equation}

This labeling system, inspired by the Wolfram code \cite{wolfram2002new}, is such that a rule number in its binary form displays the behavior of the rule. Starting from the right, the 8 first digits indicate 
the future state for each configuration as they have been ordered previously. The 8 following digits show when a division occurs. 

\begin{figure}[htbp]
\centerline{\includegraphics[width=11cm]{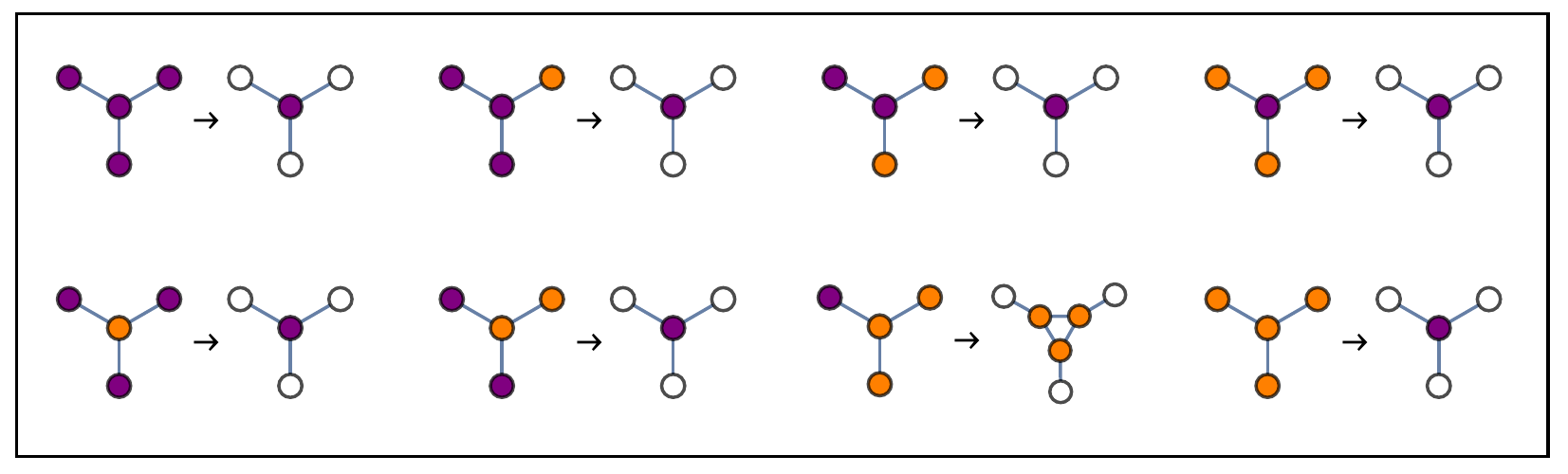}}
\caption{Rule of the example in \protect\textit{Figure \protect\ref{example}} ($n=765=10,1111,1101_2$).}
\label{example-rule}
\end{figure}

\pagebreak

\section{Implementation}

Implementing GRA as described below must leverage a sparse array format to be \linebreak efficient. \\

The application of a rule to a graph will come in several steps which are strictly equivalent to the rule being applied all at once. $o$ being the order of the graph and @ being the operator applying a function to every element of a vector, the first step consists of computing a \textbf{configuration vector} $\mathcal{C}$. 
\begin{equation}
\mathcal{C}=
\begin{pmatrix} c(v_1) \\ \vdots \\ c(v_o) \end{pmatrix}
=4\times\mathcal{S}+\mathcal{A}\cdot\mathcal{S}
\end{equation}

\noindent $\mathcal{S}$ can then be updated as follows.
\begin{equation}
    \mathcal{S}=R\,\text{@}\,\mathcal{C}
\end{equation}

\noindent A \textbf{division vector} $\mathcal{D}$ is computed similarly.
\begin{equation}
    \mathcal{D}=R'\text{@}\,\mathcal{C}
\end{equation}

Divisions can now be performed one by one as a combination of simple operations on matrices. The first 1 in $\mathcal{D}$ signals the vertex to divide. For the state vector, suffice to triple the line corresponding to this vertex. For the adjacency matrix, both the line and the column must be tripled. Ones then have to be spread across these lines and columns and the intersection must be filled with a sub-matrix containing zeros on the diagonal and ones everywhere else. Finally, the first 1 in the division vector is turned into a 0 and then tripled. This process must be repeated until $\mathcal{D}$ is a null vector. \\

\begin{figure}[H]
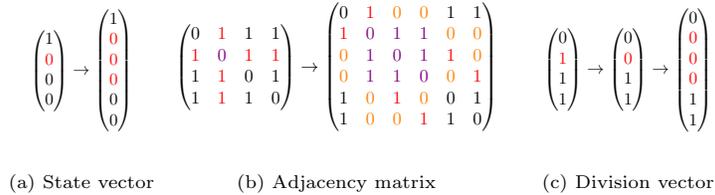

\begin{center}
\begin{subfigure}[b]{0.16\textwidth}
\begin{center}
\resizebox{0.7\textwidth}{!}{$
\begin{pmatrix}
    1 \\ {\color{red}0} \\ 0 \\ 0 
\end{pmatrix}
\rightarrow
\begin{pmatrix}
    1 \\ {\color{red}0} \\ {\color{red}0} \\ {\color{red}0} \\ 0 \\ 0
\end{pmatrix}$}
\end{center}
\caption{State vector}
\end{subfigure}
\begin{subfigure}[b]{0.38\textwidth}
\begin{center}
\resizebox{0.93\textwidth}{!}{$
\begin{pmatrix}
    0 & {\color{red}1} & 1 & 1 \\
    {\color{red}1} & {\color{violet}0} & {\color{red}1} & {\color{red}1} \\
    1 & {\color{red}1} & 0 & 1 \\
    1 & {\color{red}1} & 1 & 0
\end{pmatrix}
\rightarrow
\begin{pmatrix}
    0 & {\color{red}1} & {\color{orange}0} & {\color{orange}0} & 1 & 1 \\
    {\color{red}1} & {\color{violet}0} & {\color{violet}1} & {\color{violet}1} & {\color{orange}0} & {\color{orange}0} \\
    {\color{orange}0} & {\color{violet}1} & {\color{violet}0} & {\color{violet}1} & {\color{red}1} & {\color{orange}0} \\
    {\color{orange}0} & {\color{violet}1} & {\color{violet}1} & {\color{violet}0} & {\color{orange}0} & {\color{red}1} \\
    1 & {\color{orange}0} & {\color{red}1} & {\color{orange}0} & 0 & 1 \\
    1 & {\color{orange}0} & {\color{orange}0} & {\color{red}1} & 1 & 0
\end{pmatrix}$}
\end{center}
\caption{Adjacency matrix}
\end{subfigure}
\begin{subfigure}[b]{0.24\textwidth}
\begin{center}
\resizebox{0.77\textwidth}{!}{$
\begin{pmatrix}
    0 \\ {\color{red}1} \\ 1 \\ 1 
\end{pmatrix}
\rightarrow
\begin{pmatrix}
    0 \\ {\color{red}0} \\ 1 \\ 1 
\end{pmatrix}
\rightarrow
\begin{pmatrix}
    0 \\ {\color{red}0} \\ {\color{red}0} \\ {\color{red}0} \\ 1 \\ 1
\end{pmatrix}$}
\end{center}
\caption{Division vector}
\end{subfigure}
\end{center}
\caption{Division operated on the second vertex of the graph from \protect\textit{Figure \protect\ref{example-graph}}.}
\end{figure}

\section{Exploration}
\label{section:exploration}

An exhaustive exploration will be performed on the set of rules with exactly one division. This set being color symmetric, if it is the case as well for $G_0$, considering the half of it were the division is applied to a dead cell will give the full picture. This makes a total of 1024 rules to be explored.
\begin{equation}
    n\in\{i+2^j\;|\;i\in [[0,255]],\;j\in [[8,11]]\}
\end{equation}

\begin{figure}[H]
\centerline{\includegraphics[width=.3\textwidth]{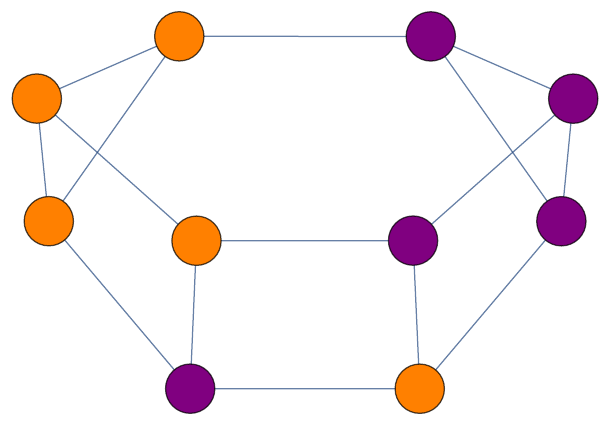}}
\caption{$G_0$: Initial graph with all 8 configurations.}
\label{initial-graph}
\end{figure}

To investigate the behavior of these rules, a simple graph $G_0$ depicted in \textit{Figure \ref{initial-graph}} \linebreak and containing all 8 configurations will serve as a starting point. This ensures that $G_1\neq G_0$.

All the results and plots in this paper were computed in Mathematica 13.0 \cite{Mathematica}, on a MacBook Pro 2020 with an Intel Core i5 Quad-Core 2.0GHz and 16GB of memory. Graphs are displayed by the default graph plotting function in Mathematica and can be better seen by enlarging this pdf.\\

\section{Growth}

One natural way to classify rules is by observing the pattern of growth they produce. \textit{Table \ref{growth-classification}} shows a classification obtained using several methods, mainly the least squares method.

\begin{table}[H]
\caption{Number Of Rules By Growth Pattern.}
\centering
\resizebox{11cm}{!}{
\begin{tabular}{|c|c|c|c|c|c|c|c|c|}
\hline
\textbf{Halted}&\multicolumn{3}{|c|}{\textbf{Linear}}&\textbf{Quadratic} &\textbf{Exponential} &\textbf{Unclassified} \\
\cline{2-4} 
 &\textbf{strict}&\textbf{periodic}&\textbf{chaotic}& & & \\
\hline
422 & 73 & 19 & 3 & 1 & 374 & 132 \\
\hline
\end{tabular}}
\label{growth-classification}
\end{table}

\subsection{Halted growth}

A total of 422 rules lead to a halt when reaching a cycle. These cycles have periods of 1, 2, 3, 6, or 8.
\begin{table}[H]
\caption{Period distribution in rules leading to cycles.}
\centering
\resizebox{6.5cm}{!}{
\begin{tabular}{|c|c|c|c|c|c|}
\hline
\textbf{Period} & 1 & 2 & 3 & 6 & 8\\
\hline
\textbf{Number of rules} & 310 & 102 & 2 & 6 & 2 \\
\hline
\end{tabular}}
\label{period-distribution}
\end{table}

\subsection{Linear growth}

73 rules produce a strictly linear growth, with the same number of vertices gained at each time step. 19 rules create a repeating pattern of growth that averages to linearity. The largest period spanned 42 time steps and was produced by rule 2183. This comes from the two main structures of the graph: a loop growing in 7 steps $\{2,0,0,2,0,0,0\}$ and a braid growing in 6 steps $\{4,0,0,4,2,0\}$.
\begin{figure}[H]
\begin{center}
\begin{tabular}{c c} \parbox{12cm}{\includegraphics[width=12cm]{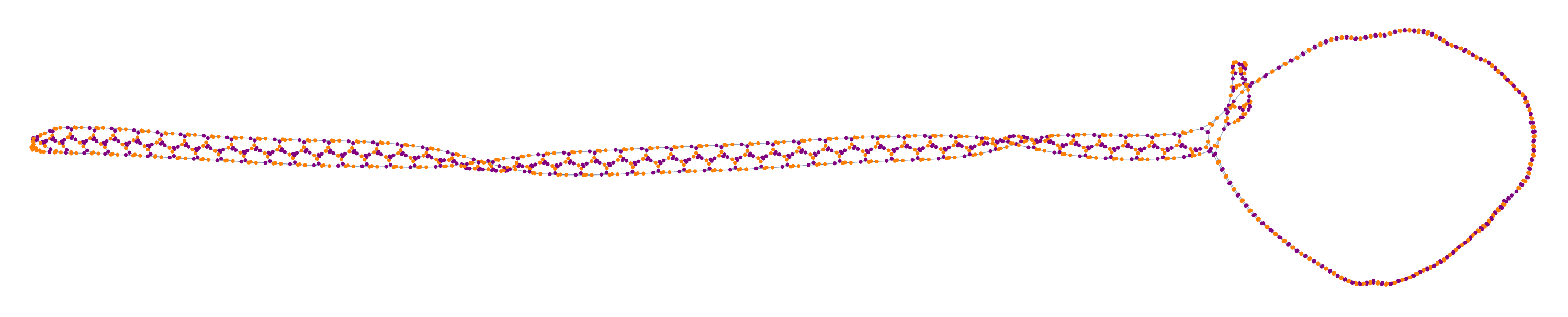}} 
 \qquad\qquad\quad \\ 
 \parbox{10cm}{\includegraphics[width=10cm]{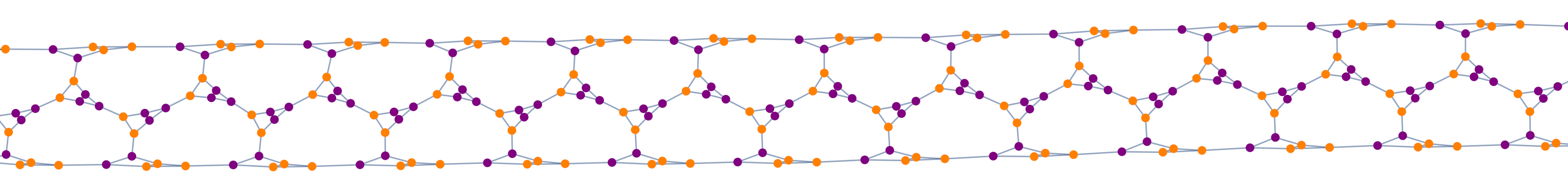}} 
 \qquad\qquad\qquad\qquad\qquad\qquad \\
\end{tabular}
\end{center}
\caption{Rule 2183 at $t=600$ along with a close-up of the braid.}
\label{2183-structure}
\end{figure}

More surprising though are the 3 rules giving rise to a chaotic linear growth: 2222, 2238 and 2239. These rules produce ladder-like loops in which states change according to an Elementary Cellular Automaton rule \cite{weisstein2002elementary}. For the sake of concision, only rule 2222 will be presented. \\

\begin{figure}[H]
\centerline{\includegraphics[width=.55\textwidth]{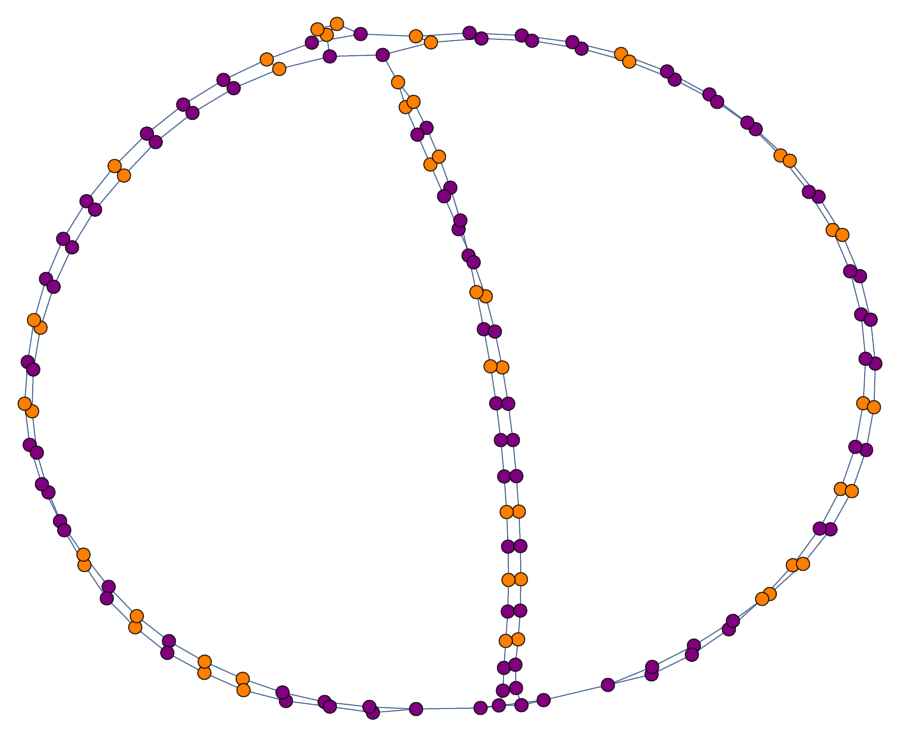}}
\caption{Rule 2222 at $t=200$.}
\label{2222-structure}
\end{figure}

This rule, although producing a rather simple graph structure, is interesting by the rhythm of its growth. At any time step, it will either grow of 0, 2, 4 or 6 new vertices which averages to a linear growth. It can be noted that there are no increase of 6 vertices at once for $t\in[[1217,\,20608]]$.
\begin{figure}[H]
\centerline{\includegraphics[width=.55\textwidth]{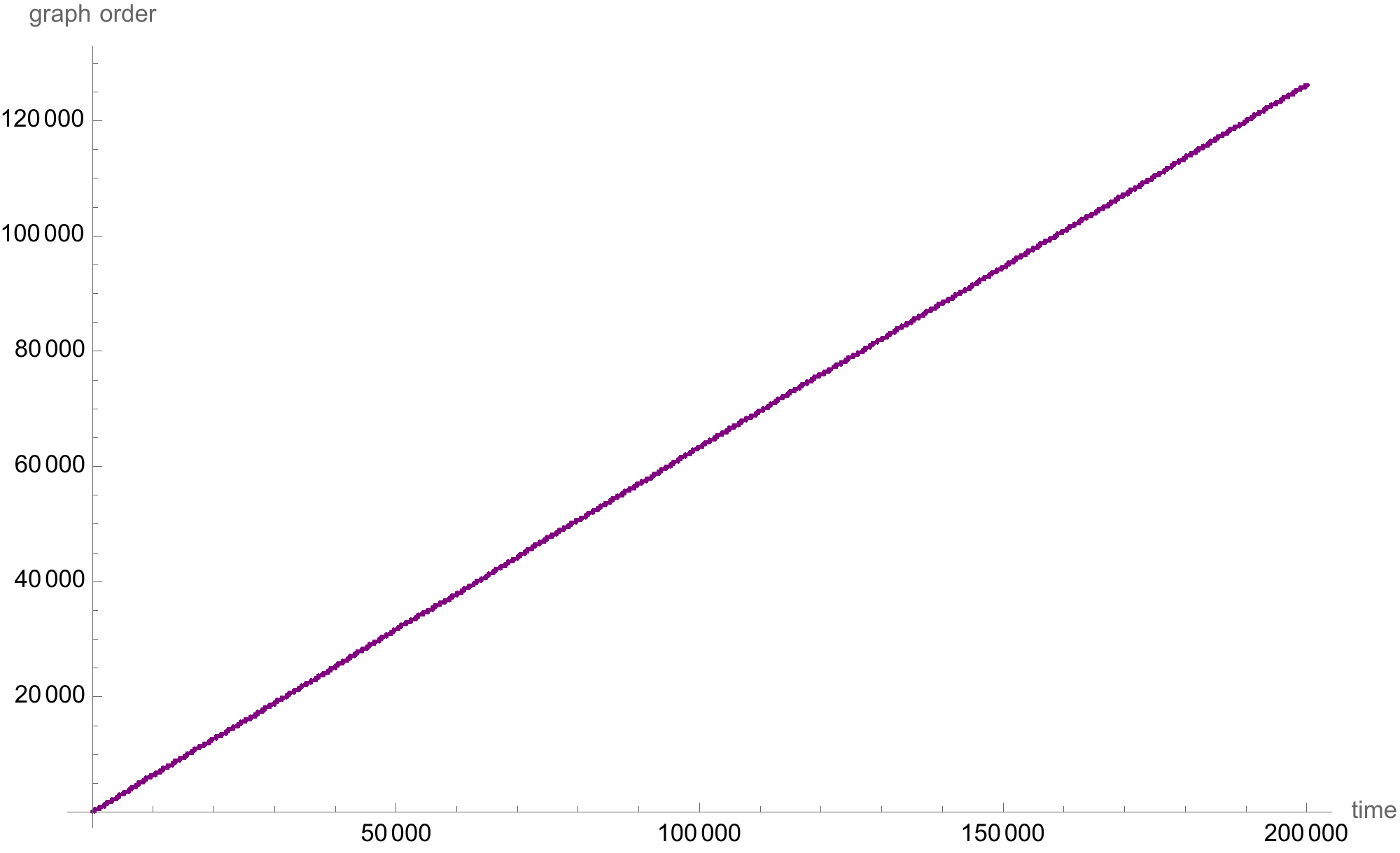}}
\caption{Rule 2222: linear growth.}
\label{2222-growth}
\end{figure}
\begin{center}
\boxed{\text{Model fit: $o = 72.8095 + 0.630684\times t$ 
\qquad\qquad Adjusted $R^2$: $0.999992$}}
\end{center}

\begin{figure}[H]
\centerline{\includegraphics[width=.57\textwidth]{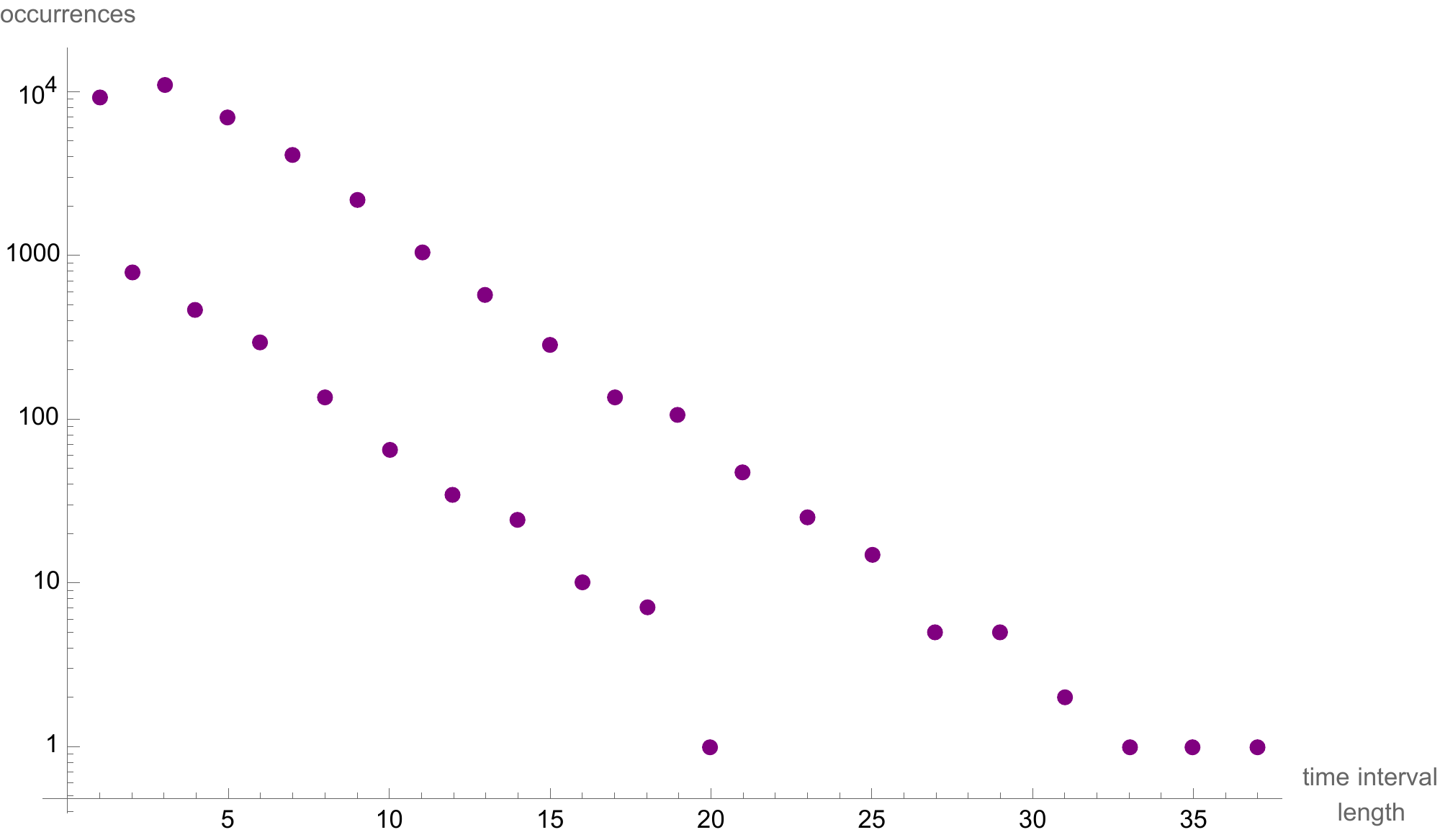}}
\caption{Rule 2222: length distribution of time intervals without growth\protect\hfill \\  computed for $0\leq t\leq 200,000$.}
\label{2222-intervals}
\end{figure}

What is remarkable here is the fact that this growth does not follow a repeating pattern but, instead, seems to have chaotic properties. To illuminate the peculiar aspect of this phenomenon, \textit{Figure~\ref{2222-intervals}} shows the distribution of the lengths of time intervals without growth in a log-linear plot.

\subsection{Quadratic growth}
\begin{figure}[H]
\centering
\begin{tabular}{c c c} 
 \parbox{1.5cm}{\includegraphics[width=1.5cm,angle=180]{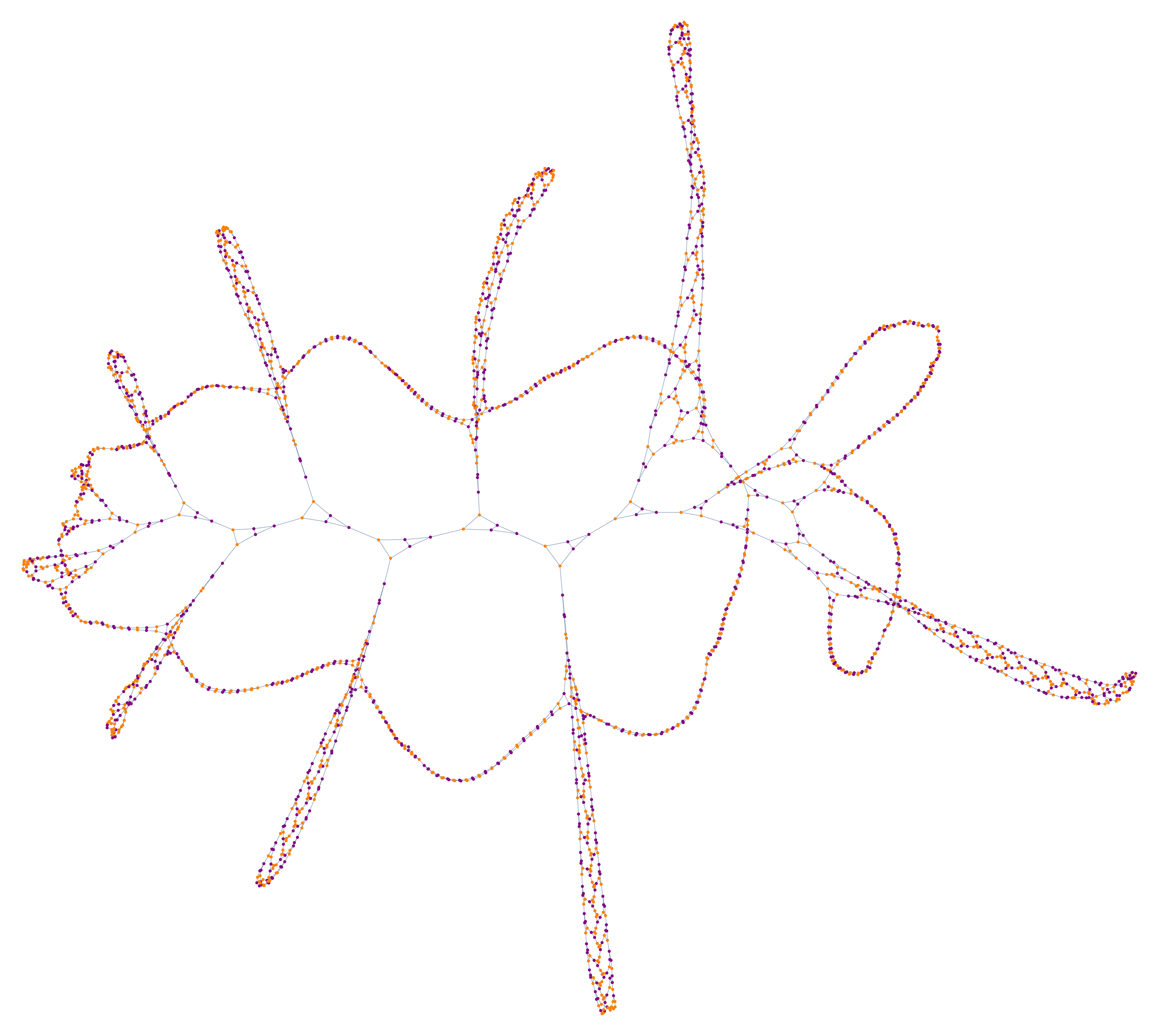}}
 \qquad & 
 \parbox{3cm}{\includegraphics[width=3cm,angle=180]{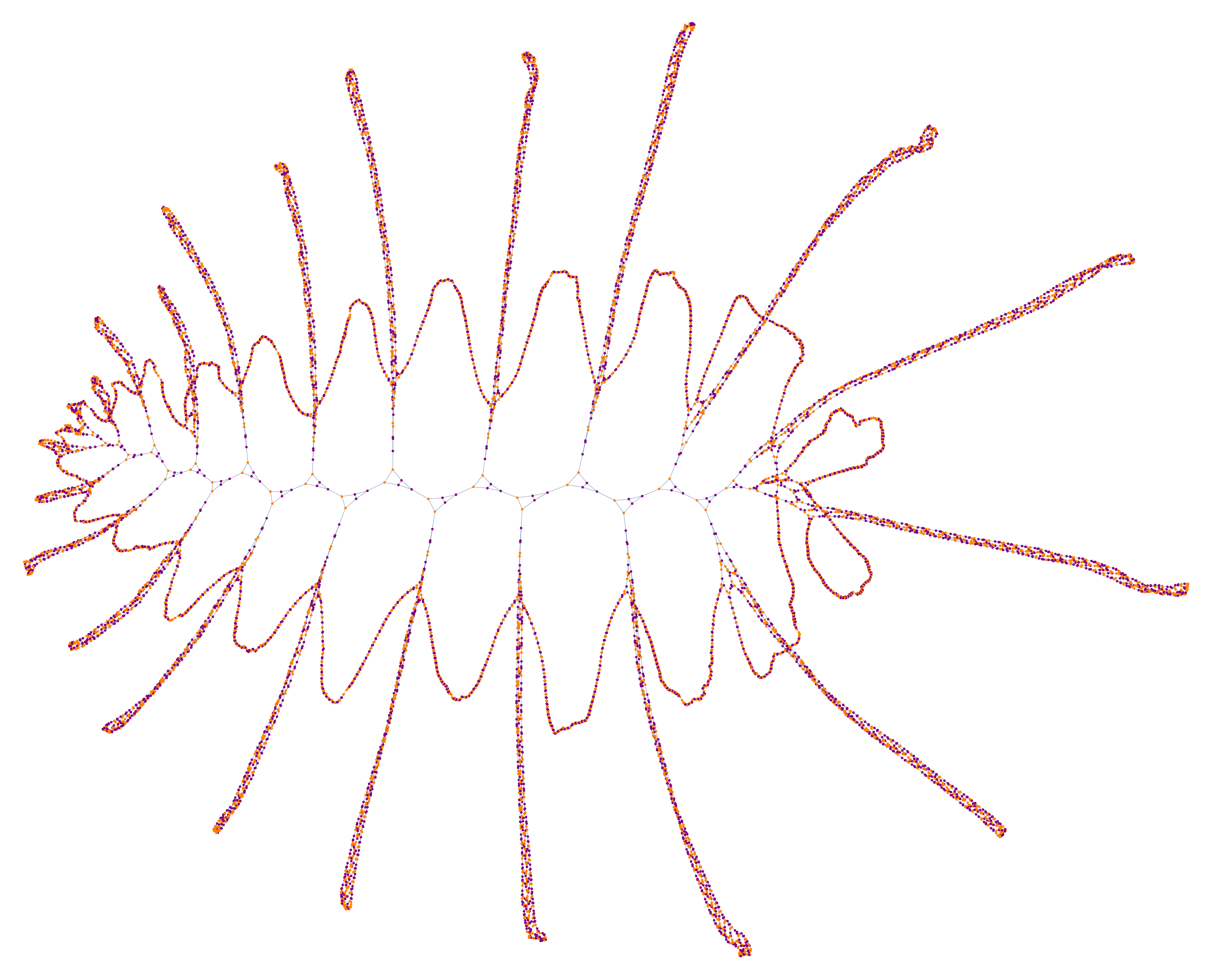}} & 
 \parbox{6cm}{\includegraphics[width=6cm,angle=180]{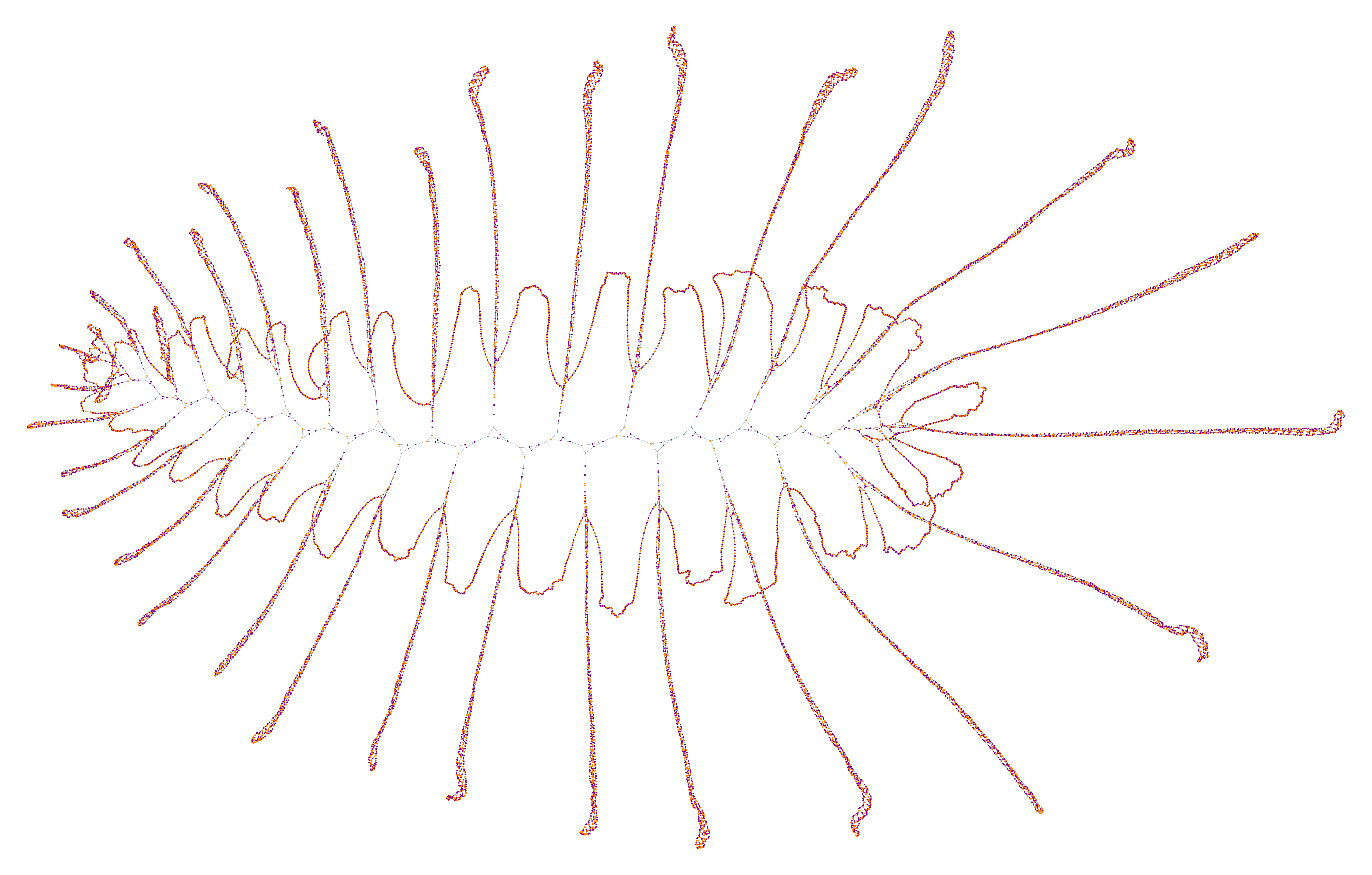}} \\
 $t=200$&$t=400$&$t=600$
\end{tabular}
\caption{Rule 2182.}
\label{2182-structure}
\end{figure}

Rule 2182 was found to be unique in the set regarding its growth. It produces \linebreak a particular structure that grows almost quadratically, the exponent being slightly above 2.

\begin{figure}[H]
\centerline{\includegraphics[width=.7\textwidth]{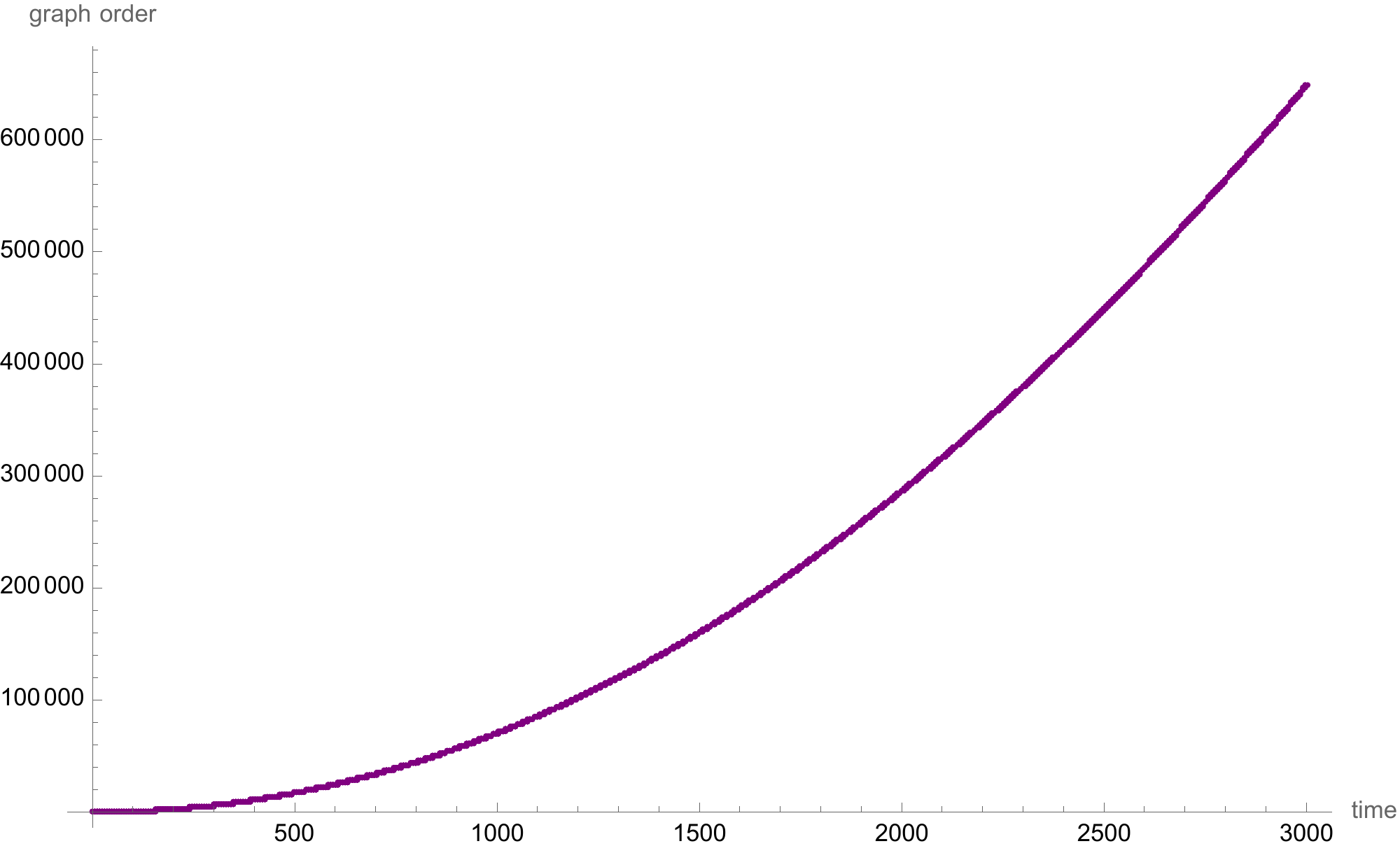}}
\caption{Rule 2182: quasi-quadratic growth.}
\label{2182-growth}
\end{figure}
\begin{center}
\boxed{\text{Model fit: $o = 0.0620341\times t^{2.01874}$
\qquad\qquad Adjusted $R^2$: $0.999999$}}
\end{center}

\vspace{1.5mm}

\subsection{Exponential growth}

Rule 256 produces a perfect exponential growth from $t=1$. Among the 373 other exponential rules, some display simple fractal structures, others are more complex.

\begin{figure}[H]
\centering
\begin{tabular}{c c}
 \parbox{60mm}{\includegraphics[height=40mm,angle=0]{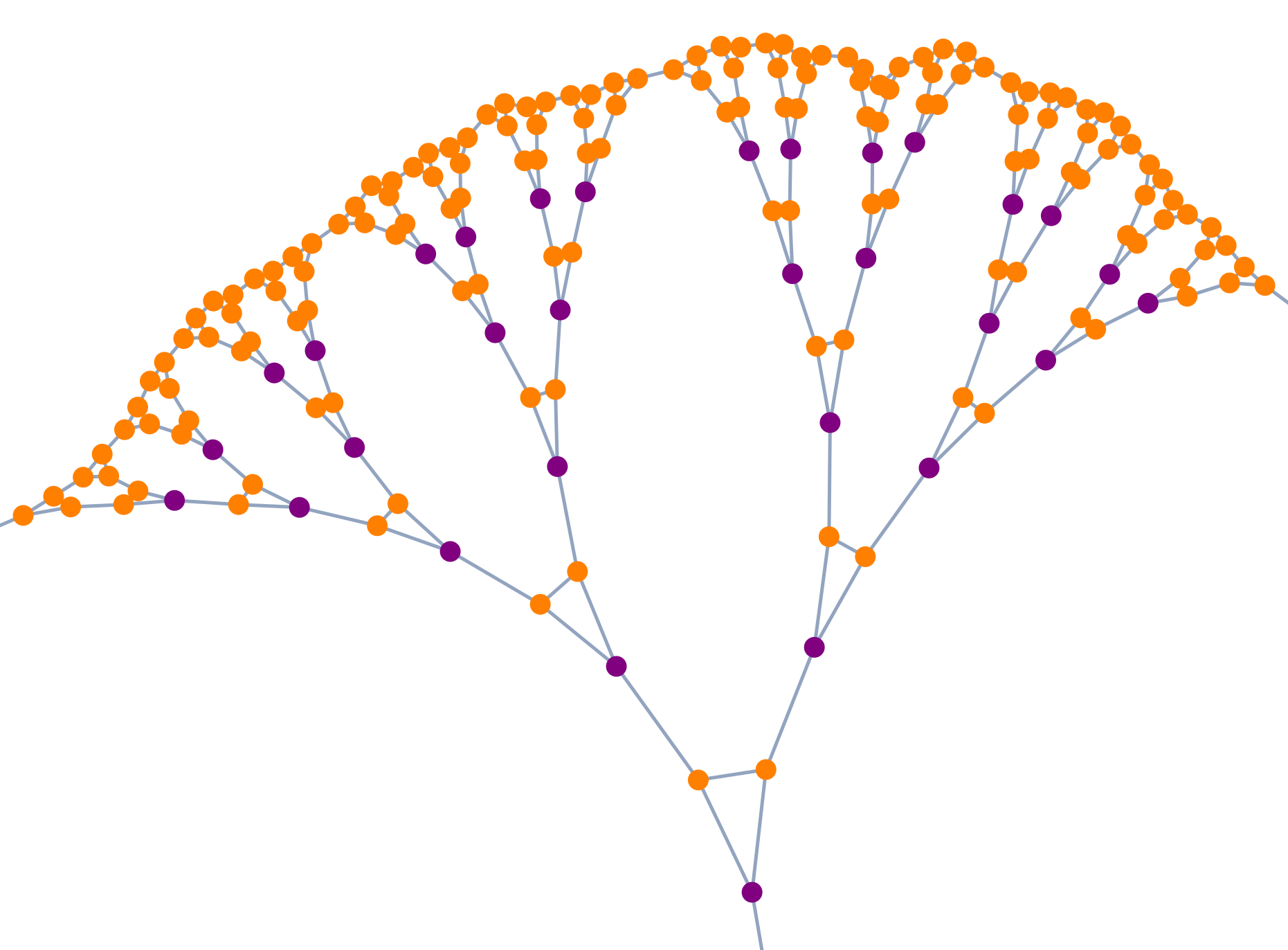}} & 
 \parbox{60mm}{\includegraphics[height=50mm,angle=90]{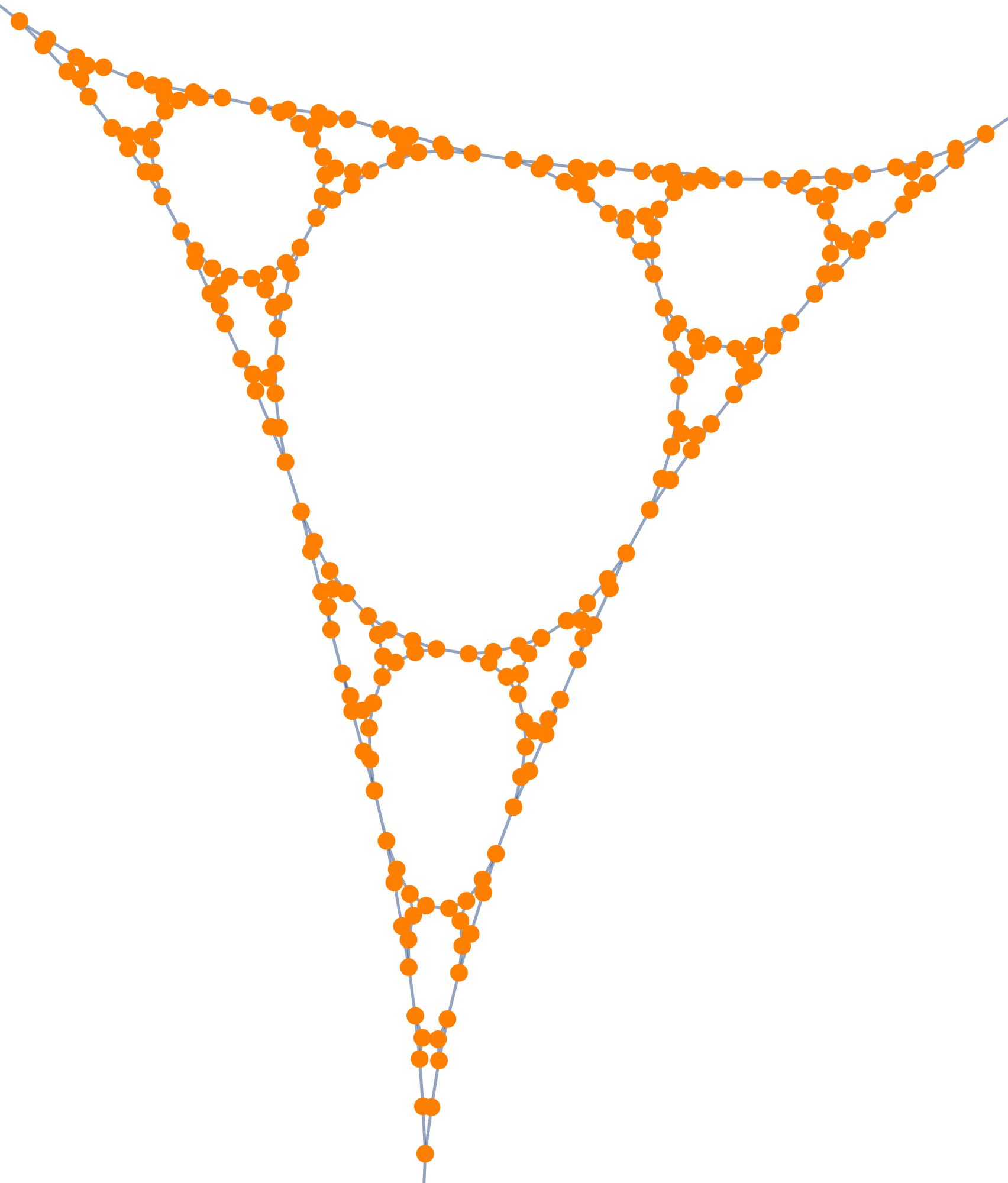}}
\end{tabular}
\caption{Examples of fractal structures found in exponential rules.}
\label{fractal}
\end{figure}

\vspace{1.5mm}

\subsection{Unclassified}

132 rules have not allowed for a confident classification. Suspecting interesting behaviors to be responsible for these less trivial growth patterns, all of them were observed individually and a selection of rules producing remarkable graph structures is presented in \textit{Section~\ref{section:structures}}. \\

\pagebreak

\section{Graph structures}
\label{section:structures}

\begin{figure}[H]
    \centering
    \begin{tabular}{c c}
         \includegraphics[width=.46\textwidth]{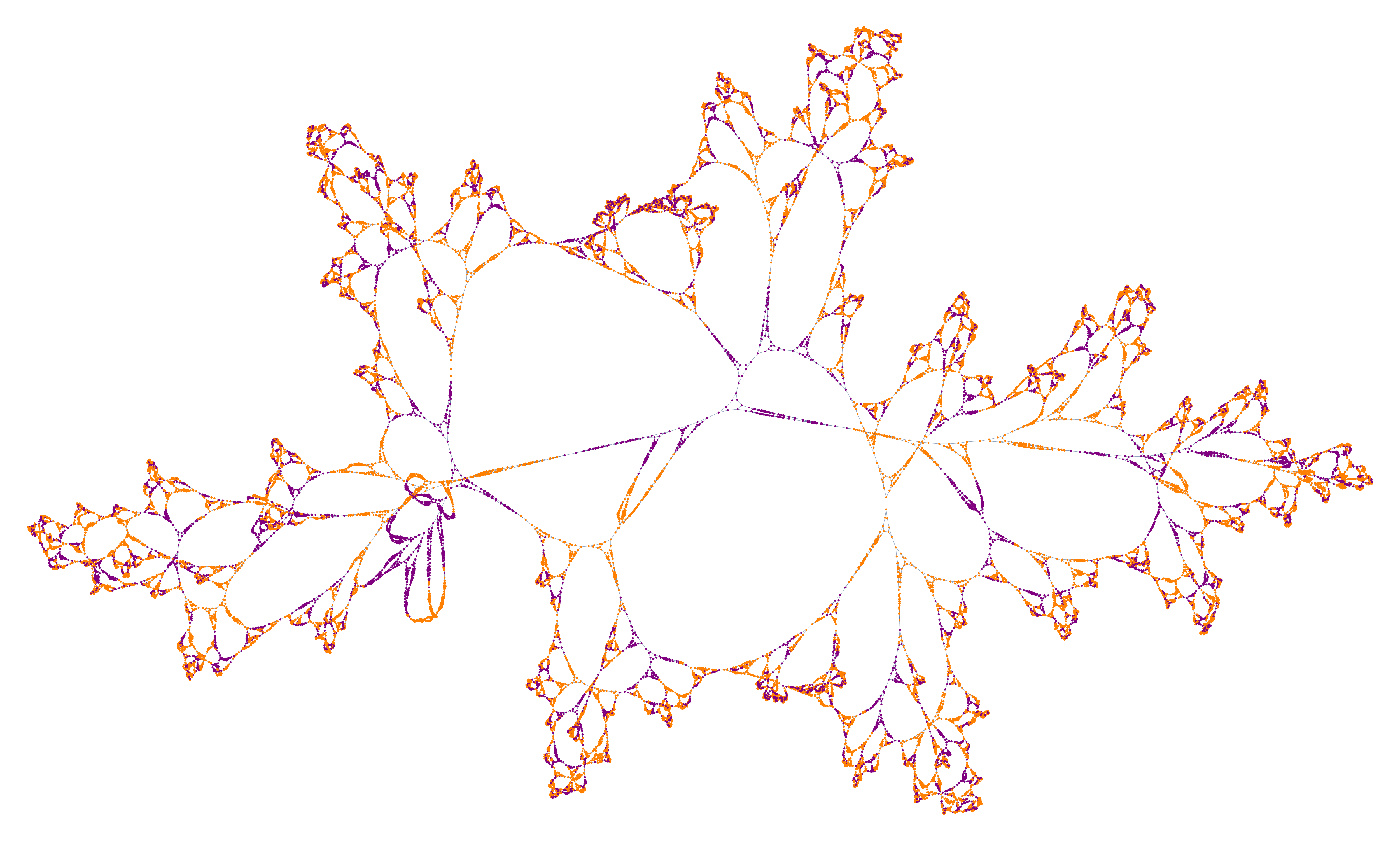} & \includegraphics[width=.46\textwidth]{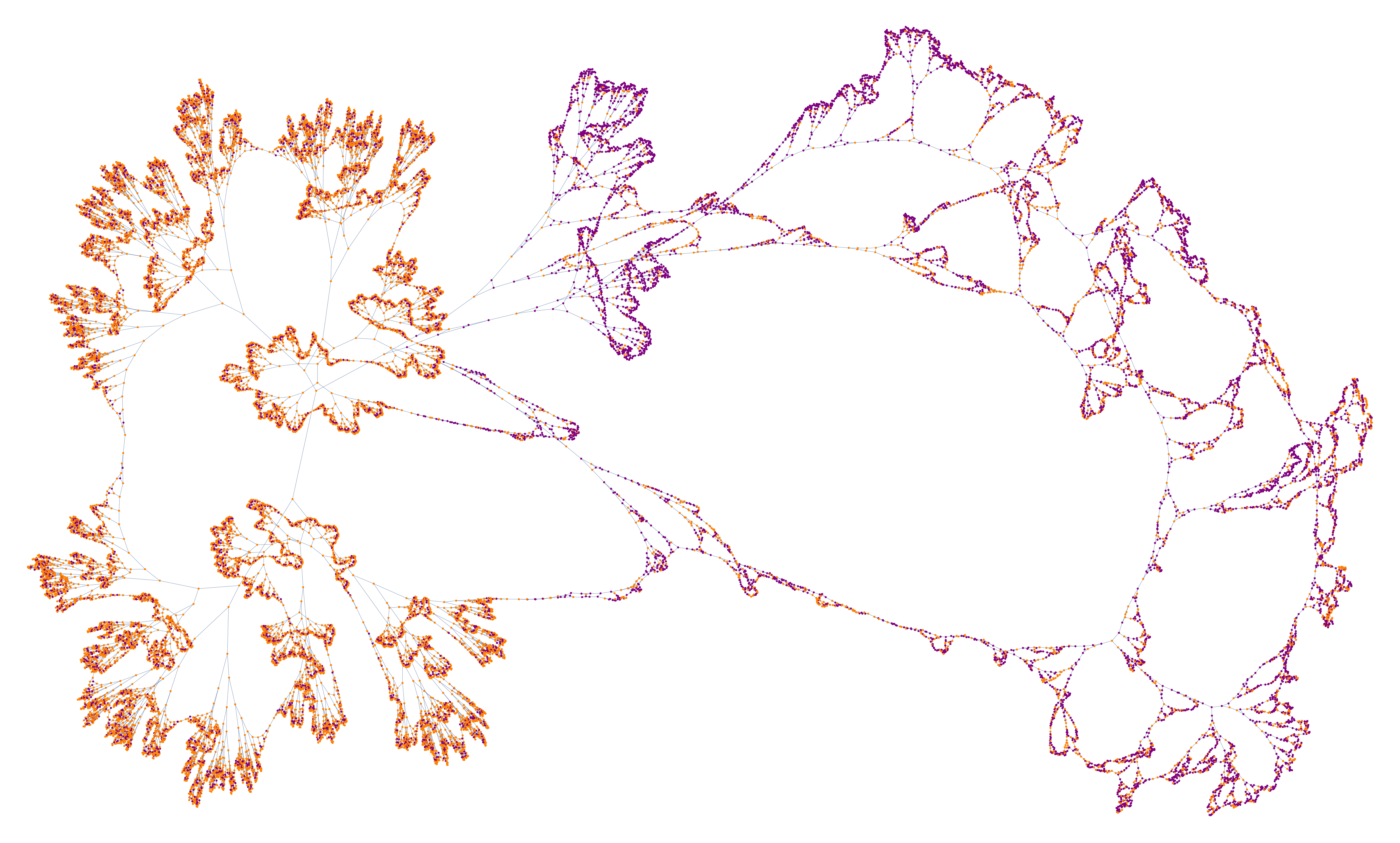} \\
         549 & 618 \\
         \includegraphics[width=.46\textwidth]{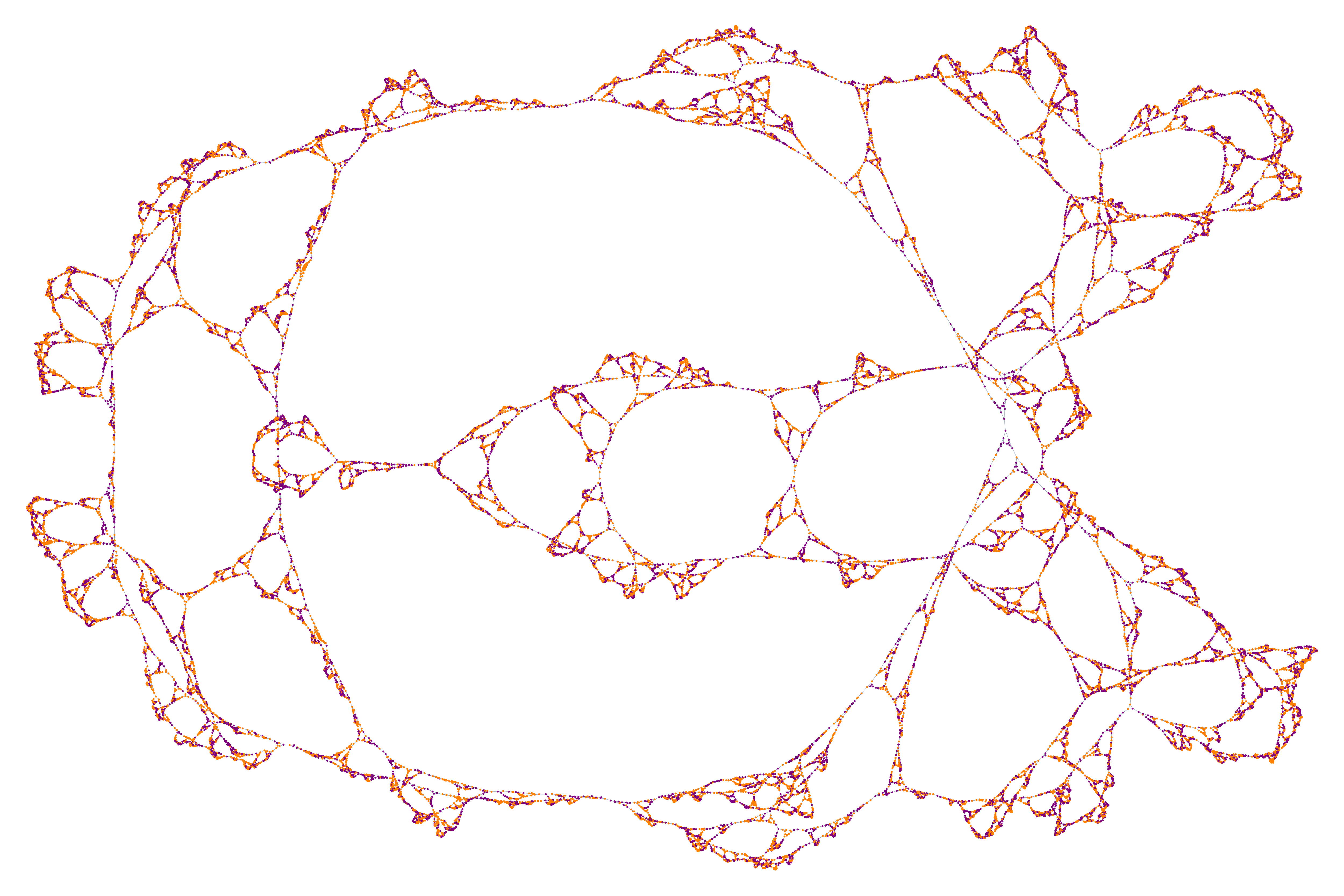} & \includegraphics[width=.40\textwidth]{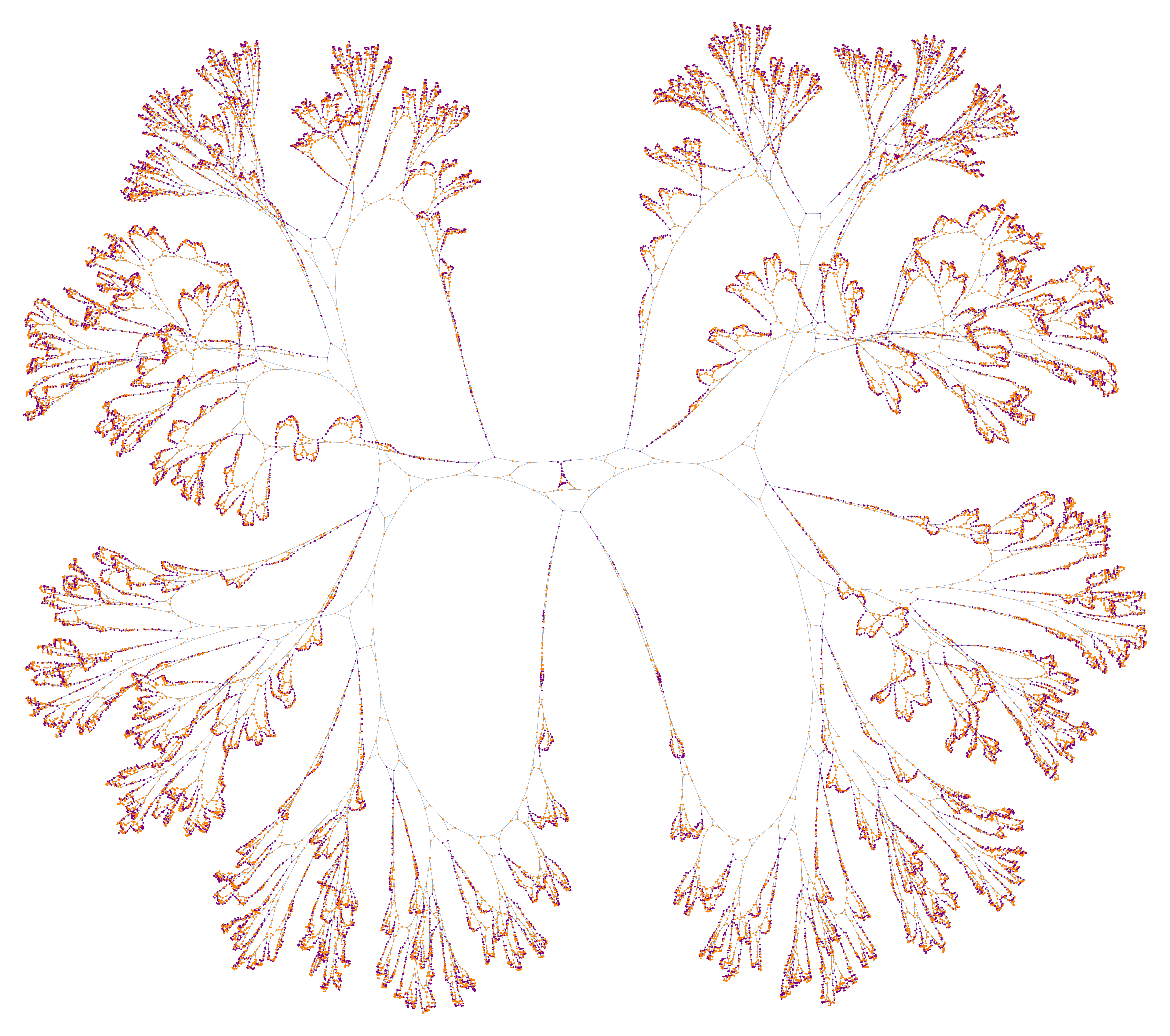} \\
        1062 & 1111 \\
        \includegraphics[width=.46\textwidth]{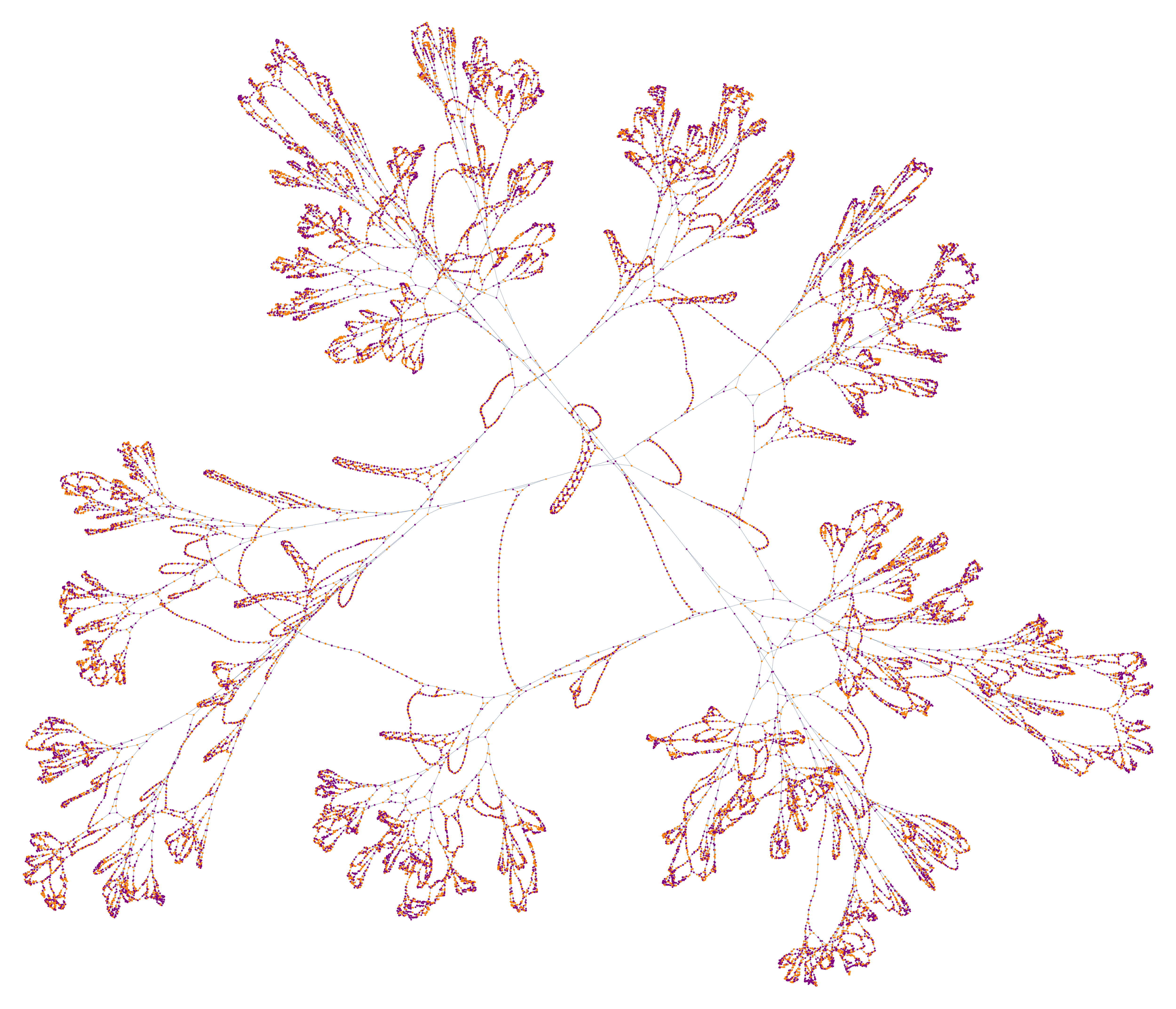} & \includegraphics[width=.46\textwidth]{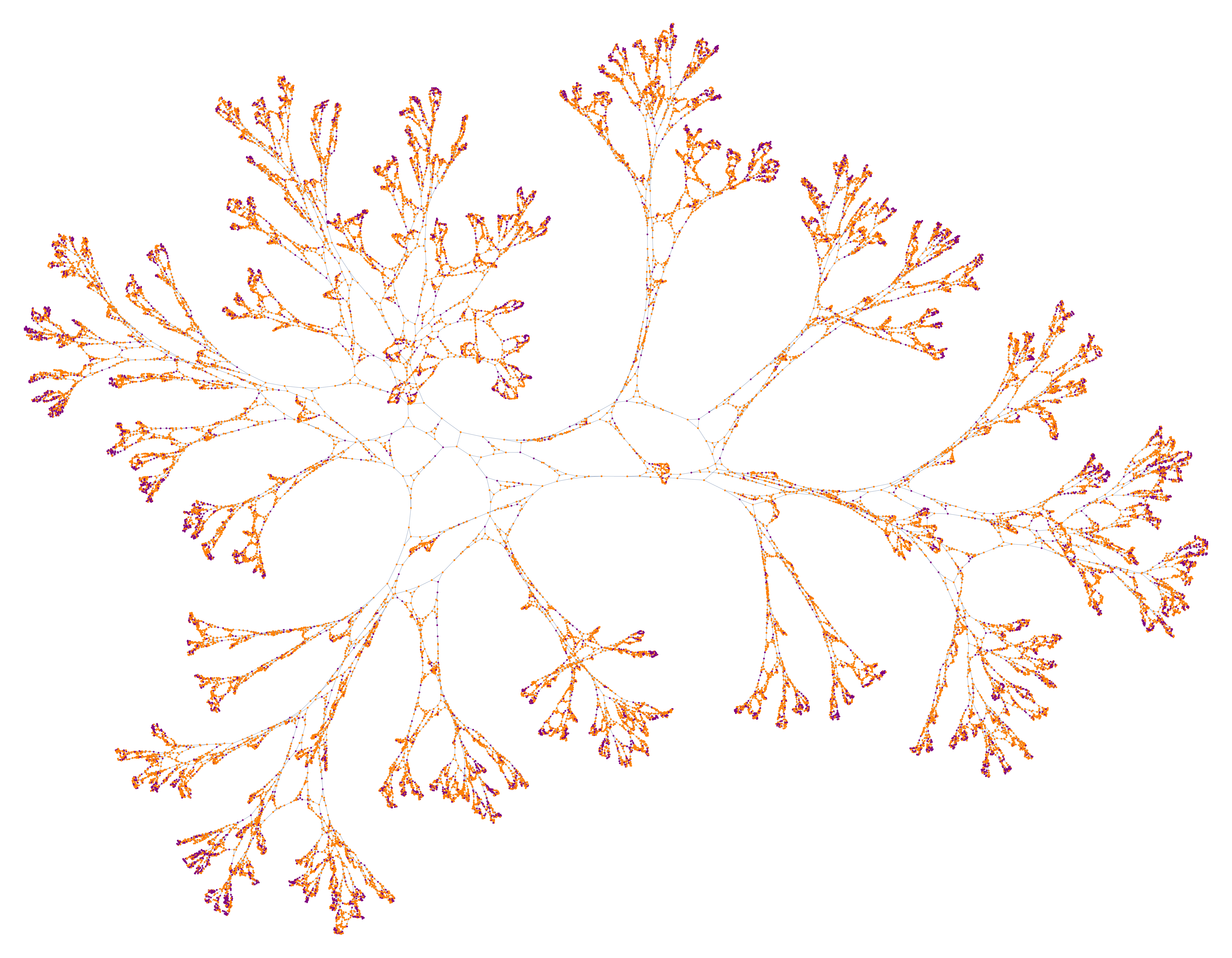}\\
        2199 & 2236
    \end{tabular}
\caption{Some organic-looking graph structures.}
\label{organic-structures}
\end{figure}

As \textit{Figure \ref{2182-structure}} already hinted at, some rules create organic-looking graph structures. \textit{Figure \ref{organic-structures}} shows a selection of graphs with this property, labeled with the number of the rule that produced them. These structures are particularly evocative. Some even seem familiar, resembling macroscopic algae or related photosynthetic eukaryotes. \\

\pagebreak

\section{Conclusion}
In this paper, a new approach to Graph-Rewriting Automata based on linear algebra was presented. When leveraging a sparse array format, it can be used to implement GRA much more efficiently than with previous methods. To demonstrate the benefits of this approach, a natural subset of GRA was explored only using a MacBook Pro 2020. \\

A focus was first put on the evolution of the number of vertices over time, leading to a classification of rules by growth pattern. Several noteworthy phenomena were discovered through this approach. It also led to a category of unclassified growths in which interesting behaviors were suspected. The graphs of this category were thus visually inspected on a individual basis and some remarkable structures were discovered. \\

In future work, the ideas presented here could easily be adapted to cover a variety of other types of GRA. Here are a few examples of possible trajectories:
\begin{enumerate}
    \item[•] using continuous valued or non-binary discrete states,
    \item[•] working with d-regular or non-regular graphs,
    \item[•] reaching an extended neighborhood with powers of $\mathcal{A}$,
    \item[•] working with directed graphs using $\mathcal{A}\cdot\mathcal{S}$ and $\mathcal{A}^T\cdot\mathcal{S}$,
    \item[•] including other topology altering operations.
\end{enumerate}

\noindent Other tools to analyse graphs directly from their adjacency matrix and state vector could also be developed to gain new insights into the behaviors and properties of GRA. \\

\section*{Acknowledgment}
This work was partially supported by the Université Grenoble Alpes, through an \linebreak Excellence Internship which took place between June 7\textsuperscript{th} and July 15\textsuperscript{th} 2022. \\

\bibliographystyle{IEEEtran}
\bibliography{refs.bib}

\end{document}